\documentclass[reqno]{amsart}

\usepackage[all]{xy}

\numberwithin{equation}{section}

\theoremstyle{plain}
\newtheorem{lemma}{Lemma}[section]
\newtheorem{teo}[lemma]{Theorem}
\newtheorem{propo}[lemma]{Proposition}
\newtheorem{coro}[lemma]{Corollary}

\newtheorem{claim}{Claim}

\theoremstyle{definition}

\theoremstyle{remark}

\newcommand{\pic}{{\rm Pic\thinspace}}
\newcommand{\bs}{{\rm Bs\thinspace}}

\newcommand{\p}{\mathbb{P}}

\newcommand{\oc}{{\mathcal O}}

\newcommand{\ls}{{\mathcal L}}

\newcommand{\ts}{\tilde{S}}

\newcommand{\rt}{\longrightarrow}

\begin{document}

\title{Special systems through double points on an algebraic surface}
\author{Antonio Laface}
\address{
Antonio Laface \newline
Department of Mathematics and Statistics \newline
Queen's University \newline
K7L 3N6 Kingston, Ontario (Canada)
}
\email{alaface@mast.queensu.ca}
\keywords{Linear systems, double points, secant varieties} \subjclass{14C20}
\begin{abstract}
Let $S$ be a smooth algebraic surface satisfying the following property: $H^i(\oc_S(C))=0$ ($i=1,2$) for any irreducible and reduced curve $C\subset S$. The aim of this paper is to provide a characterization of special linear systems on $S$ which are singular along a set of double points in general position. As an application, the dimension of such systems is evaluated in case $S$ is an Abelian, an Enriques, a $K3$ or an anticanonical rational surface.
\end{abstract}
\maketitle

\section*{Introduction}
In what follows $S$ will be a smooth algebraic surface defined over an algebraically closed field of characteristic $0$.

The  problem of characterizing special linear systems $|H|$ through $k+1$ double points in general position on $S$ is strictly connected with the problem of evaluating the dimension of the $k$-secant variety of $S$. That is why this subject and its generalizations have been studied by many authors (see for example~\cite{cgg,cc,cm1,da,te}). The main results on this subject are related to the classification of defective surfaces, i.e. surfaces whose $k$-secant variety is defective. This means that $H$ is assumed to be very ample and even in that case it is not easy to understand which are the numerical characters of the special pair $(S,H)$. Trying to fill this gap, this paper is mainly devoted to the study of linear  
systems on those surfaces which share the following property:
\begin{equation}\label{h^i}
H^i(\oc_S(C))=0\hspace{5mm}{\rm for}\hspace{5mm} i=1,2
\end{equation}
for any irreducible and reduced curve $C\subset S$. 
The  main result is Proposition~\ref{num}. As an application a complete characterization of special linear systems on Abelian, Enriques, $K3$ and anticanonical rational surfaces is given. \\
The paper is divided as follows: 
In section 1 some preliminary material about linear systems and $k$-secant varieties 
is  given together with a partial classification of surfaces satisfying~\eqref{h^i}.
Section 2 deals with the main part of the paper, where it is stated and proved 
the characterization of these special systems. As an application, in section 3, special linear systems on Abelian, Enriques and $K3$ surfaces are completely classified. As a consequence it is proved that no one of these surfaces can have a defective $k$-secant variety. Finally section 4 focuses on the proof of Gimigliano-Harbourne-Hirschowitz conjecture for linear systems of the form $\ls_2(d;m_1,\ldots,m_9,2^r)$. As a consequence the complete list of defective anticanonical surfaces is given.

\section{Preliminaries and notation}

In what follows $S$ will be a smooth algebraic surface defines over an algebraically closed field of characteristic $0$. We will denote by $H$ an effective and integral divisor of $S$. The dimension of the cohomology groups $H^i(S,\oc_S(H))$ will be denoted by $h^i(\oc_S(H))$. 
As usual $K_S$ will denote the class of the canonical divisor of $S$.
The map $\pi: \ts\rt S$ will denote the blowing-up of $S$ along points $p_0,\ldots,p_k$ in general position and the exceptional divisors will be denoted by $E_0,\ldots,E_k$.
A linear system $|L|$ on $\ts$ will be called special if
\[
h^0(\oc_{\ts}(L))h^1(\oc_{\ts}(L)) > 0.
\]







\begin{propo}\label{class}
Let $S$ be a smooth algebraic surface for which property~\eqref{h^i} is satisfied, then $p_g(S)=q(S)=0$ unless $S$ is one of the following:
\begin{enumerate}
\item a $K3$ surface
\item an abelian surface which does not contain elliptic curves
\item an hyperelliptic surface.
\end{enumerate}
\end{propo}
\begin{proof}
If $p_g(S) = 0$ and $q(S)\geq 2$, then (see~\cite{be}) any minimal model $S'$ of $S$ is a  ruled surface. Consider the composition $\psi: S\stackrel{\pi}{\rt} S'\stackrel{\phi}{\rt} B$, where $B$ is the base of  the ruling of $S'$ and let $F=\psi^{-1}(p)$ for a general $p\in B$.
From the fundamental sequence of $F$
\[
\xymatrix{ 0 \ar[r] & \oc_{S}\ar[r] & \oc_{S}(F) \ar[r] & \oc_{F}(F) \ar[r] & 0}
\]
one obtains that $h^1(\oc_S(F)) = q(S)-1 > 0$.
If $q(S)=1$ and $S'$ is ruled, then $\phi$ has a section $C'$ with $C'^2\leq 0$. Taking the fundamental sequence of $C'$ one obtains that $h^1(\oc_{S'}(C'))>0$ and this in turn implies that $h^1(\oc_S(C))>0$ where $C$ is the strict transform of $C'$. 
\\
Observe that if $K_S-C$ is effective for some integral curve then by Serre duality $h^2(\oc_S(C)) = h^0(\oc_S(K_S-C)) > 0$, so that $K_S$ is effective only if  it is trivial. 
In this case $S$  is either a $K3$  or an abelian surface. In the second case it is easy to show that if $S$ contains an elliptic curve $C$ then either $h^1(\oc_S(C))>0$ or $|2C|$ contains an integral curve with the same property.
\end{proof}

\section{Linear systems through double points on surfaces}\label{msec}

In this section we will consider the blow-up $\ts$ of a smooth algebraic surface $S$ along 
$k$ general points. With $H$ we will denote the pull-back of an integral curve of $S$.
We will denote by $L$ the divisor
\[
H-\sum_{i=0}^{k-1} 2E_i,
\]
where $E_i = \pi^{-1}(p_i)$ are the exceptional divisors of the blow-up. 
We will assume that $|L|$ is a non-special system and that $|L-2p|$ is special for 
the general $p\in S$. 
Let 
\[
|L| = F + |M|
\]
be the decomposition of $|L|$ into its fixed and component-free part.

\begin{propo}\label{curve}
The system $|M|$ is composed with a linear pencil $|D|$ and $M=nD$ with $n>1$.
\end{propo}
\begin{proof}
First of all observe that if $|L-2p|$ is special for the general $p\in\ts$, then the same is true for the system $|M-2p|$ since $\dim |L| = \dim |M|$ and $p$ can be chosen outside the base locus of $|L|$. 
Now let $C={\rm im\thinspace}\phi_{|M|}(S)$ and let $T_q(C)$ be the tangent space at a smooth point of $C$, then, since $T_q(C)$ imposes at most $2$ conditions on the hyperplanes containing it, $\dim T_q(C) =1$ and this means that $C$ is a curve.
Observe that if $C$ is non-rational then $\phi$ must be a morphism, since otherwise one of the rational curves obtained resolving the indeterminacy of $\phi$ would dominate $C$. After lifting $\phi$ to the normalization $\tilde{C}$ of $C$ and applying the Stein factorization to this lift, we obtain a morphism $\ts\rt C'$ with connected fibres. Since $g(C')>0$ then $q(\ts)>0$ which, by Proposition~\ref{class}, means that $\ts$ is the blowing-up of either a $K3$ or an hyperelliptic surface. In both these cases $\chi(\oc_S)=0$ so that by Riemann-Roch theorem applied to $\oc_{\ts}(D)$ one has that $h^0(\oc_{\ts}(D)) = -\frac{K_{\ts}D}{2} < 0$, which is absurd.\\
Bertini's second theorem (see~\cite{kl}) implies that that either the general element of $|M|$ is irreducible or it is composed with a pencil. In the first case $\dim |M| =1$ and $|M-2p|$ would be empty for $p$ general. This implies that $|M|=|nD|$ with $n>1$.
\end{proof}

In this way it is easy to observe that given $|L-2p|$ must contain a double curve through $p$. That is why Proposition~\ref{curve} implies this well-known corollary:

\begin{coro}\label{double}
A special linear system through double points has a double curve as a fixed component.
\end{coro}

Let us start investigating the properties of the base locus of $|L|$.

\begin{lemma}\label{bs}
No one of the $E_i$ is contained into the base locus of $|L|$.
\end{lemma}
\begin{proof}
Let us consider the divisor $R = L + 2E_i$, then by hypothesis $|R-2E_i| = |R-rE_i|$ with $r\geq 3$. Since $E_i$ is the blowing-up of a general point of $S$, this implies that 
\[
|R-2p| = |R-rp|
\]
for the general $p\in S$. Now, we want to prove that the preceding relation implies that the system $|R-2p|$ is empty, which means that $|L|$ is empty. Let $|R| = R_{\rm fix} + |R_{\rm free}|$ be the decomposition of tthe system into its fixed and free part.
By the generality assumption on $p$, we can always assume that $p\not\in R_{\rm fix}$ and this means that the system $|R_{\rm free}|$ must have the  same property of $|R|$.
So we can always assume that $|R|$ has no fixed locus. 
Let $V\subseteq\p^N$ be the image of the map $\phi_{|R|}$. \\
If $\dim V = 1$ then the assumption on $R$ is equivalent to ask that $V$ has infinitely many flexes and this is possible only if $V$ is a line. But in this case $|R-2p|$ is empty. \\
If $\dim V = 2$ then, as before, one can deduce that the general hyperplane section of $V$ is a curve with infinite flexes and do it is a line. But the only surface containing a family of dimension $N\geq 2$ of lines is the plane and this still implies that $|R-2p|$ is empty.
\end{proof}

\begin{lemma}\label{02}
For any exceptional divisor $E_i$ the intersection $M E_i$ is $0$ or $2$. 
\end{lemma}
\begin{proof}
From the decomposition $|L| = F + |M|$, one obtains that 
\[
2 = E_i F + E_i M.
\]
By lemma~\ref{bs} $E_i$ is not contained into $\bs |L|$, then $E_iF \geq0$ and this gives $0\leq E_i M\leq 2$. Now observe that if $E_iM=1$ then $1=ME_i=nDE_i$ would give $n=1$ and this is not possible by proposition~\ref{curve}.
\end{proof}


\begin{lemma}\label{vanish}
Let $S$ be a smooth algebraic surface for which property~\eqref{h^i} is satisfied and such that $ME_i = 2$ for at least one $i$, then
one has that $h^1(\oc_{\ts}(rD)) = 0$ for $r=1,2$.
\end{lemma}
\begin{proof}
The divisor $D\in\pic(\ts)$ can be written as $H'-\sum E_i$, where $H'$  is the pull-back of a divisor of $S$. Since $D$ is irreducible and reduced, then the same is true for the general element of $|H'|$. By Hypothesis $h^1(\oc_S(H'))=0$ and since $|D|$ corresponds to a linear system on $S$ through simple points in general position, then its dimension is the expected one and this means that $h^1(\oc_{\ts}(D))=0$. \\
About $|2D| = |M|$, observe that it is fixed component free, hence by corollary~\ref{double} it is non-special.
\end{proof}

The preceding lemma allows one to find the numerical characters of the curve $D$ by means of Riemann-Roch theorem. 

\begin{propo}\label{num}
Let $S$ be a smooth algebraic surface for which property~\eqref{h^i} is satisfied and such that $ME_i = 2$ for at least one $i$, then
the general element of $|D|$ is a smooth curve and
\[
D^2 =  \chi(\oc_S)-1,\hspace{1cm} DK_{\ts} = 3\chi(\oc_S)-5.
\]
Moreover, for  any irreducible and reduced $C\subseteq\bs(L)$ one has that $h^1(\oc_{\ts}(C))=0$ and $CD = D^2$. 
\end{propo}
\begin{proof}
The equalities $h^0(\oc_{\ts}(D)) = 2$ and $h^0(\oc_{\ts}(2D)) = 3$, together with the vanishing of the higher cohomology groups, imply that
\[
\frac{D^2-DK_{\ts}}{2} + \chi(\oc_S) =  2\hspace{1cm}{\rm and}\hspace{1 cm}
\frac{4D^2-2DK_{\ts}}{2} + \chi(\oc_S) =  3.
\]
Solving these equations for $D^2$ and $DK_{\ts}$ one obtains the numerical properties of $D$. By Proposition~\ref{class} it follows that $\chi(\oc_S)\leq 2$ which in turn implies that $D^2\leq 1$ so that $|D|$ can have at most one simple base point $p$. By Bertini's first theorem the general element of $|D|$ is smooth away from $p$ and, obviously, it has to be smooth also in $p$, since otherwise two elements of $|D|$ would have a bigger intersection at that point. \\
The curve $C$ can be written as $H'' -\sum E_i -\sum 2E_j$ for some indexes $i$ and $j$, where $H''\in\pi^*\pic(S)$.  Since $C$ is irreducible and reduced the same is true for the general element of $H''$, so that by assumption $h^1(\oc_{\ts}(H''))=0$. The vanishing of $h^1(\oc_{\ts}(C))$ follows now from Corollary~\ref{double}. \\
From the exact sequence 
\[
\xymatrix{ 0 \ar[r] & \oc_{\ts}\ar[r] & \oc_{\ts}(D) \ar[r] & \oc_{D}(D) \ar[r] & 0}
\]
and $D^2\leq 1$ one obtains that $h^0(\oc_D(D)) = 1$ and $q(S)=0$. Tensoring the preceding sequence with $\oc_{\ts}(C)$ and using the vanishing of $h^i(\oc_{\ts}(C))$ for $i=1,2$, one obtains that $h^0(\oc_D(C+D)) = 1$ and $h^1(\oc_D(C+D)) = 0$. This by Riemann-Roch implies that 
\[
1 = (C+D)D+1-g(D),
\]
which is equivalent to say that $CD = D^2$.
\end{proof}


\begin{propo}\label{zero}
Let $S$ be a smooth algebraic surface for which property~\eqref{h^i} is satisfied and such that $ME_i = 0$ for any $i$ then $n | 2 (1-\chi(\oc_{S}))$. Moreover if $\chi(\oc_S)=1$ then $D^2=0$ and the general element of $|D|$ is a smooth curve.
\end{propo}
\begin{proof}
The hypothesis implies that $D\in\pi^*\pic(S)$ so that by assumption $h^i(\oc_{\ts}(D)) = 0$ for $i=1,2$. This together with the relation $\dim |nD| = n$ gives
$(n^2D^2-nDK_{\ts})/2+\chi(\oc_{\ts}) = n+1$, which implies the first part of the thesis. \\
If $\chi(\oc_S)=1$ then the preceding equality gives $nD^2-DK_{\ts} = 2$ which, together with the equality, $D^2-DK_{\ts}=2$ implies that $D^2=0$. This imply that $|D|$ is base point free and by Bertini's first theorem its general element is smooth.
\end{proof}

\begin{coro}\label{p=q=0}
Let $S$ be a smooth algebraic surface with $p_g(S)=q(S)=0$ for which property~\eqref{h^i} is satisfied, then $D^2=0$ and $S$ is rational.
\end{coro}
\begin{proof}
The fact that $D^2=0$ is an easy consequence of Propositions~\ref{num} and~\ref{zero}. In both cases $h^i(\oc_{\ts}(D))=0$ for  $i=1,2$ and this means that $D^2-DK_{\ts}=2$. By the genus formula the general element of $|D|$ is a smooth rational curve, hence $S$ is rational.
\end{proof}

\section{Applications to some non-rational surfaces}

In this section we apply what has be done to some surfaces of Kodaira dimension $0$. In what follows $S$ will be either a $K3$ or an abelian or an Enriques surface. The following is the main property that will be used in what follows:

\begin{lemma}\label{k0}
Let $H$ be an irreducible and reduced divisor on $S$, then $h^1(\oc_{\ts}(H)) = 0$ unless $S$ is an abelian surface and $H$ is an elliptic curve or $S$ is an Enriques surface and $p_a(H)\leq 1$.
\end{lemma}
\begin{proof}
If $H^2 > 0$ then, since $H\sim K_S + (K_S + H)$ and $K_S + H$ is nef and big, the result is achieved by applying the Kawamata-Vieveg vanishing theorem.
The other possibilities are $H^2=0,-2$ and in both cases the fundamental sequence of 
the general element of $|H|$ gives the thesis.
\end{proof}

The preceding proposition allows one to apply results from section~\ref{msec} to these surfaces.
The first theorem prove the conjecture stated in~\cite{dl} for linear systems through double points.

\begin{teo}\label{teo-k3}
Let $S$ be a $K3$ surface, then the system  $|L-2E_k|$ is special if and only if $k=1$ and $L=2H-2E_0$ with $H^2=2$.
\end{teo}
\begin{proof}
From Proposition~\ref{num} one has that $D^2=1$ and $DK_{\ts}=1$.
Observe that by Proposition~\ref{zero} together with Proposition~\ref{curve} implies that $n=2$ even if $ME_i = 0$ for any $i$. This means that 
$D$ can always be written as $D = H-\sum_{i=1}^r E_i$ where $H\in\pi^*\pic(S)$.
The canonical divisor $K_{\ts} = \sum_{i=0}^kE_i$ since $K_S=\oc_S$. 
From these two equalities one obtains
\[
1 = DK_{\ts} = (H-\sum_{i=1}^r E_i)(\sum_{i=1}^r E_i) = r,
\]
which gives $r=1$. From $D^2=1$ one deduce that $H^2 = 2$. \\
Now let us suppose that $|L|$  has a fixed part and let $C\subseteq\bs(L)$ be an irreducible and reduced curve. Observe that if $D = H - E_0$, then $C E_0 = 0$  since $2 = L E_0 = F E_0 + 2(D E_0)$. Hence $C$ can be written as 
$C = R-\sum_{i=1}^s a_i E_i$,
where $R\in\pi^*\pic(S)$ and  $a_i=1,2$.
By Proposition~\ref{num} $CD=1$ and so also $RH=1$. This means that $\phi_{|H|}(R)$ is a line, where $\phi_{|H|}$ is the double cover of $\p^2$, and that $R$ is rational.
So $R$ is the pull-back of a $(-2)$-curve on $S$ and this implies that $s=0$ by the generality assumption on the $E_i$'s. In order to conclude,
observe that if $C\in\pic(\ts)$ is a $(-2)$-curve then $2D+C_{|C} \cong \oc_{\p^1}$, so that the vanishing of $h^1(\oc_{\ts}(2D))$ implies that $\dim |2D+C| = \dim |2D| + 1$ and $C\not\subseteq\bs(L)$.
\end{proof}

Observe that in this case $2H$ is not very ample since the map $\phi_{|2H|}$ factorizes through $\phi_{H}$ and the $2$-Veronese embedding of $\p^2$ (that is where the speciality comes from).



Let us consider the case of an abelian surface:

\begin{teo}\label{abelian}
Let $S$ be an abelian surface, then the system  $|L-2E_k|$ is never special. 
\end{teo}
\begin{proof}
Let $D$ be as in Proposition~\ref{curve}, then by Proposition~\ref{num}
has that, if $h^1(\oc_{\ts}(D))=0$ then $D^2=-1$, which is not possible since $\dim |D| = 1$. 
Then the only possibility if that $h^1(\oc_{\ts}(D))>0$, which by lemma~\ref{vanish} implies that
$D\in\pic(\ts)$ can be written as $H'-\sum E_i$, where $H'$  is the pull-back of a divisor of $S$ and $h^1(\oc_{\ts}(H'))>0$. By lemma~\ref{k0} $H'$ is an elliptic curve and this implies that
$\dim |H'| =0$ which is a contradiction.
\end{proof}

\begin{teo}\label{enriques}
Let $S$ be an Enriques surface, then the system  $|L-2E_k|$ is never special.
\end{teo}
\begin{proof}
Following the first lines of the proof of Proposition~\ref{num}, one see that $h^1(\oc_{\ts}(D))>0$, since otherwise $|D|$ would be a pencil of rational curves.
Writing $D = H'-\sum E_i$ and observing that $h^1(\oc_{\ts}(H'))=h^1(\oc_{\ts}(D))$, by lemma~\ref{k0} we conclude that $|H'|$ is an elliptic pencil and that $D=H'$. 
Now let $C$ be an integral curve contained in the fixed part of $|L|$, we proceed by assuming the following:
\begin{claim}\label{bscurve}
$CE_i>0$ for at least one $E_i$
\end{claim}
Write $C=H''-\sum c_iE_i$ where $1\leq c_i\leq 2$ and observe that if $CH'=0$ then  $H''H'=0$ which implies that $H''=H'$ and that $C=H'-E_i$.
On the other hand if $CH'>0$ then $C$ does not live in a fiber of $|H'|$ so, either $h^1(\oc_{\ts}(C)) = 0$ or $C$ is of the form $H''-E_i$ for some elliptic pencil $|H''|\neq |H'|$ ($C$ can not be a $(-2)$-curve by the general position assumption on the points). In the first case, looking at the exact sequence
\[
\xymatrix{ 0 \ar[r] & \oc_{\ts}(C)\ar[r] & \oc_{\ts}(C+H') \ar[r] & \oc_{H'}(C+H') \ar[r] & 0}
\]
a dimension count shows that $CH'\leq 1$. On the other hand the preceding is not possible since $H'\in\pic(S)$ is an even class. By a similar argument one can argue that in the second case $CH'\geq 4$, which gives a contradiction since $h^1(\oc_{\ts}(C)) = h^1(\oc_{\ts}(H'')) = 1$.
\end{proof}

\begin{proof}[Proof of Claim~\ref{bscurve}]
Suppose that $C$ does not intersect any $E_i$, then it is the pull-back of a divisor of $S$. By Riemann-Roch theorem one has that $C^2\leq 0$, so that $C$ can be either a $(-2)$-curve or an elliptic curve with $h^1(\oc_{\ts}(C))=0$. Since $C+2H'$ is nef and big then, by Kawamata-Viewheg vanishing theorem, one has that $h^0(\oc_{\ts}(C+2H'))=4$. On the other hand $h^0(\oc_{\ts}(2H))=3$ implies that $C$ can not be in the fixed part of $|L|$. 
\end{proof}

\section{Applications to some anticanonical rational surfaces}
In this section $S$ will be the blowing-up of $\p^2$ along at most $9$ points in general position. We will denote by $e_i$ the exceptional divisors of the blowing-up of the nine points and by $h$ the pull-back of a general line of $\p^2$.
As proved in~\cite{har1}, an irreducible and reduced curve $H\in\pic(S)$ satisfies the conditions~\eqref{h^i}.
The main theorem of this section is a proof of the Gimiggliano-Harbourne-Hirschowitz conjecture 
for a class of linear systems:

\begin{teo}\label{plane}
A linear system $|L-2E_k|$ of the form 
\[
|dh - \sum_{i=1}^9 m_i e_i -\sum_{i=0}^k 2E_i|
\]
is  special if and only if there exists a $(-1)$-curve $E$ such that $(L-2E_k)R\leq -2$.
\end{teo}
\begin{proof}
One part of the proof is just an application of the Riemann-Roch theorem. For the ``only if'' part, let $H= dh - \sum_{i=1}^9 m_i e_i$ and observe that if $|H|$ is special then, as proved in~\cite{har1}, there exists a $(-1)$-curve $E$ of type $\delta h-\sum_{i=1}^9 \mu_i e_i$ such that $(L-2E_k)E = HE \leq -2$. \\
If $|L|$ is non-special and $|L-2E_k|$ is, then writing $|L| = F + |M|$, by Proposition~\ref{num}, one has that $M=2D$ with $D^2=0$ and $DK_{\ts} = -2$. This Implies that the curve $R:=D-E_k$ is a $(-1)$-curve. \\
Now , let $C\subseteq F$ be an irreducible and reduced curve. Still, by Proposition~\ref{num}
one knows that $CD = 0$, which gives 
\[
LD = (F+2D)D = 0.
\]
This allows one to calculate the intersection
\[
(L-2E_k)R = (L-2E_k)(D-E_k) = -2,
\]
which proves the thesis.
\end{proof}

As an application of Theorem~\ref{plane}, one can evaluate the dimension of the secant variety of $S\subset\p^n$. First of all we briefly recall from~\cite{har2} a condition that will be required in what follows. A divisor $H=dh-\sum_{i=1}^9 m_ie_i$ is in {\em standard form} if $m_1\geq\ldots\geq m_r\geq 0$ and $d\geq m_1 + m_2 + m_3$.

\begin{teo}\label{sec-p2}
Let  $H=dh-\sum_{i=1}^9 m_ie_i$ be a very ample divisor on the blowing-up $S$ of the plane along nine points in general position, then the $k$-secant variety of $\phi_{|H|}(S)$ is defective if and only if $H$ is one of the following:
\[
2h\hspace{5mm}4h\hspace{5mm}2nh-(2n-2)e_1
\]
with $k=1, 4$ and $2n-1$ respectively.
\end{teo}

Before proving the theorem we need to prove some lemmas. 
Let $\tilde{\p}^2$ be the blowing-up of $\p^2$ along $r$ points in general position and let 
$C=\delta h-\sum_{i=1}^r \mu_i e_i$ be a curve on $\tilde{\p}^2$. The degree of $C$ is defines as $\deg C=\delta$.

\begin{lemma}\label{noether}
Let $C$ be a $(-1)$-curve in $\tilde{\p}^2$, then there exists a quadratic Cremona transformation $\sigma$ such that $\deg \sigma(C) < \deg C$.
\end{lemma}
\begin{proof}
A proof of this can be found in~\cite{cm2}. It depends on the fact that, if $C\in|\delta h-\sum_{i=1}^r \mu_i e_i|$, then $C^2+\mu_3CK_{\tilde{\p}^2} < 0$. This inequality is equivalent to
\[
\delta(\delta-\mu_3)+\sum_{i=1}^3\mu_i(\mu_3-\mu_i)+\sum_{i=4}^r\mu_i(\mu_3-\mu_i) < 0
\]
where the last sum of the first member is non-negative. By assuming that $\delta\geq \mu_1+\mu_2+\mu_3$ one can substitute this value for $\delta$ in the first part of the inequality obtaining a contradiction.
\end{proof}

\begin{lemma}\label{intersection}
Let $\ls=|dh-\sum_{i=1}^r m_ie_i|$ be a linear system in standard form, then $\ls C\geq 0$ for any $(-1)$-curve $C$.
\end{lemma}
\begin{proof}
Let $C\in|\delta h-\sum_{i=1}^r \mu_i e_i|$ be a $(-1)$-curve and
let $\sigma$ be a quadratic transformations which decreases the degree of $C$.
Consider the intersection
\begin{eqnarray*}
\ls\sigma(C) & = & (dh-\sum_{i=1}^r m_ie_i) ((\delta-t) h-\sum_{i=1}^3 (\mu_i-t) e_i-\sum_{i=4}^r \mu_i e_i) \\
& = & d\delta - dt -\sum_{i=1}^rm_i\mu_i + t\sum_{i=1}^3m_i \\
& = & \ls C -t(d-\sum_{i=1}^3 m_i),
\end{eqnarray*}
then $\ls\sigma(C)\leq\ls C$. Still by lemma~\ref{noether} there exists a sequence
of curves $C_i=\sigma_i(C_{i-1})$ for $i=1,\ldots,n$ such that $C_0=C$ while $C_n\in\ls_2(1;1^2)$ and this implies that $\ls C\geq\ls C_n = d-m_1-m_2\geq 0$.
\end{proof}

\begin{proof}[Proof of Theorem~\ref{sec-p2}]
By hypothesis $H=dh-\sum_{i=1}^9 m_ie_i$ is a very ample divisor on a rational surface whose anticanonical divisor is effective. We can always assume that the multiplicities are ordered in this way: $m_1\geq\ldots\geq m_9$.
As proved in~\cite[Theorem 2.1]{har2}, a divisor $H$ is very ample if and only if it is in standard form and the following numerical conditions are satisfied:
\[
d\geq m_1 + m_2 + 1\hspace{5mm}{\rm and}\hspace{5mm} 3d-\sum_{i=1}^rm_i\geq 3.
\]
Since $|L-2E_k|$ is special, Theorem~\ref{plane} implies that there exists a $(-1)$-curve $R\subset\ts$ such that $(L-2E_k)R=-2$. Then by lemma~\ref{intersection} one can immediately deduce that $L-2E_k$ is not in standard form and by the preceding inequalities $m_3\leq 1$. 
Now suppose that $m_2\geq 3$, then $m_1+m_2+1\leq d < m_1+m_2+2$ implies that $L-2E_k$ is
\[
(m_1+m_2+1)h-m_1e_1-m_2e_2-\sum_{i=0}^k2E_i-\sum_{i=3}^t e_i
\]
where $2\leq t\leq 9$.
By  applying one quadratic transformation based on the first three points of $L-2E_k$, one obtains the system 
\[
(m_1+m_2)h-(m_1-1)e_1-(m_2-1)e_2-\sum_{i=1}^{r}2E_i-E_0-\sum_{i=3}^t e_i
\]
which is in standard form and this means that $L-2E_k$ can not have negative intersection with any $(-1)$-curve.
The other possibility is that $m_2\leq 2$ and in this case the linear system $|L-2E_k|$ is a quasi-homogeneous special system through double (also simple) points. These systems have been studied and classified in~\cite{cm1}. It turns out that the only special ones are:
\[
|4h-\sum_{i=1}^5 2e_i|\hspace{3mm} |2nh-(2n-2)e_1-\sum_{i=2}^r2e_i|\hspace{3mm} |dh-de_1-\sum_{i=2}^r2e_i|.
\]
The third system cannot be $|L-2E_k|$ since otherwise $|H|$ would be composed with a pencil (and this means that it can not be very ample). 
The second system gives $H=2nh-(2n-2)e_1$ and there are two possibilities for this case according to $n=1$ or $n\geq 2$. In the first case $H=2h$ and the system is very ample on $\p^2$ while in the second case the system is very ample on the blowing-up of $\p^2$ along a point. The third system gives $H=4h$.
\end{proof}

\bibliographystyle{plain}

\end{document}